\documentclass[10pt]{ijnam}
\def\currentissue{1}
\def\currentyear{2022}
\def\month{February}
\pagespan{1}{15}
\copyrightinfo{2022}{} 

\usepackage{graphicx}
\usepackage{tikz}
\usepackage{amsmath} 
\usepackage{amssymb}
\usepackage{amsthm}
\usepackage{units}   
\usepackage{float}   
\usepackage{graphicx}
\usepackage{amsfonts}
\usepackage{mathrsfs}
\usepackage[colorlinks]{hyperref}
\usepackage{framed}
\usepackage{url}
\usepackage{xcolor}
\usepackage{multirow}
\usepackage{pgfplots}
\usetikzlibrary{spy}

\newcommand{\FT}{\textsc{FronTier}}
\begin{document}

\title[Validation and Verification of Richtmyer-Meshkov Instability]
{Validation and Verification of Turbulence Mixing due to Richtmyer-Meshkov Instability of an air/SF$_6$ interface}

\author[Kaman\and Holley]{Tulin Kaman \and Ryan Holley}
\address{
  Department of Mathematical Sciences,
  University of Arkansas,
  Fayetteville, AR 72701, USA
}
\email{tkaman@uark.edu \and rh027@uark.edu}
\urladdr{https://kaman.uark.edu}



\date{February,28, 2022 and, in revised form, June 14, 2022.}


\subjclass[2000]{35Q31,76F25, 60A10}

\abstract{Turbulent mixing due to hydrodynamic instabilities
occurs in a wide range of science and engineering applications 
such as supernova explosions and inertial confinement fusion. 
The experimental, theoretical and numerical studies help us to 
understand the dynamics of hydrodynamically unstable interfaces
between fluids in these important problems. 
In this paper, we present an increasingly accurate and robust front tracking method for 
the numerical simulations of Richtmyer-Meshkov Instability (RMI)
to estimate the growth rate.
The single-mode shock tube experiments of Collins and Jacobs 2002~\cite{ColJac02} 
for two incident shock strengths ($M=1.11$ and $M=1.21$) 
are used to validate the RMI simulations. 
The simulations based on the classical fifth order 
weighted essentially non-oscillatory (WENO) scheme of 
Jiang and Shu~\cite{JiaShu96} with Yang's artificial compression~\cite{Yan88} 
are compared with Collins and Jacobs 2002 shock tube experiments. 
We investigate the resolution effects using front tracking with WENO 
schemes on the two-dimensional RMI of an air/SF$_6$ interface. 
We achieve very good agreement on the
early time interface displacement and amplitude growth rate
between simulations and experiments for Mach number $M=1.11$.
A 4\% discrepancy on early-time amplitude 
is observed between the fine grid simulation 
and the $M=1.21$ experiments of Collins and Jacobs 2002. 
}

\keywords{Turbulent Mixing, Rayleigh-Taylor Instability, 
Richtmyer-Meshkov Instability, Front Tracking,
Weighted Essentially Non-Oscillatory Scheme}

\maketitle
\section{Introduction}
Turbulent mixing due to hydrodynamic instabilities occurs in 
many scientific and engineering applications.
The formation of gravitational induced mixing in oceanography;
supernovae explosions in astrophysics, and
the performance assessment for inertial confinement fusion (ICF)
are ideal sub-problems to study and understand the dynamics of turbulence and mixing~\cite{AbaSre10}.  
In the dynamics of turbulence, hydrodynamic instabilities of fluid flows 
such as Kelvin--Helmholtz, Rayleigh--Taylor, Richtmyer--Meshkov
are observed.  The review papers of Zhou~\cite{Zho17a, Zho17b},  and
Abarzhi, Gauthier and Sreenivasan~\cite{AbaGauSre13}
provide detailed resource information on the theory, experiment and
computations of these important physical instabilities.
While the velocity difference at the interface between two fluids 
develops the Kelvin--Helmholtz instability (KHI), 
the density difference with constant and impulsive acceleration
develops the Rayleigh--Taylor instabilitiy (RTI) and
Richtmyer--Meshkov instabilitiy (RMI) respectively. 
RTI arises at the perturbed interface between two fluids of different 
densities whenever the pressure gradient opposes the density gradient. 
RMI arises when a shock wave interacts with the perturbed interface. 
RMI is also known as impulsive or shock-induced RTI. 
An overview of RTI and RMI and 
the effects of material strengths, chemical reactions and magnetic fields, 
as well as the role of the instabilities in scientific and engineering applications
can be found in two review articles \cite{ZhoClaCla19, ZhoWilRam21}.
The evolution of the perturbed interface development and the
 interaction between the fluids at the macro/meso/micro 
 length scales have been the main interest of researchers.

The motions of fluid flows are described by the Euler equations, 
a system of partial differential equations, 
where the effect of molecular processes are neglected. 
The more general equations are Navier-Stokes equations (NSE). 
These problems are deeply multiscale and the level of scales 
that are desired to be resolved identify the characteristic properties. 
The three numerical approaches to model turbulence are 
(i) Direct Numerical Simulation (DNS)~\cite{MoiMah98}, 
the full NSE are resolved without any models for turbulence, 
(ii) Large Eddy Simulation (LES), the flow field is resolved 
down to a certain length scale and scales smaller than 
that are modeled rather than resolved, and 
(iii) Reynolds-Averaged Navier-Stokes (RANS), 
the time-averaged equations are solved for mean values of all quantities. 
In the LES family,  the implicit LES (ILES) solves the governing equations using an implicit 
subgrid scale model. ILES assumes that small and unresolved scales are 
purely dissipative and the numerical discretization errors are the source of artificial dissipation~\cite{GriMarLG07}. 
These approaches are briefly reviewed in section~\ref{sec:turb}
and a recent review article summarizes in detail~\cite{Sch20}.
Within these approaches, the most accurate DNS has the 
highest computational cost and the least accurate RANS is the most
desired approach in complex engineering applications due to the low 
computational cost.  

There are hydrodynamic codes such as 
\textsl{CFDNS,
HYDRA, 
Miranda, 
RAGE,
RAPTOR,
TURMOIL,}
and \FT\
that have been under continuous development to study the 
RTI/RMI induced flow, turbulence and mixing to predict the growth rate of the mixing zone accurately
based on the mathematical and numerical frameworks~\cite{Zho17a, Zho17b}. 
\textsl{CFDNS} is a Los Alamos National Laboratory (LANL) hydrodynamic simulation code designed 
for direct numerical simulation of turbulent flows~\cite{LivMohPet09}.
\textsl{HYDRA} is based on arbitrary Lagrangian-Eulerian (ALE) mesh used 
for the numerical simulation of instabilities in ICF laser-driven hohlraum~\cite{LanKarMa14}.  
\textsl{Miranda} is a Lawrence Livermore National Laboratory (LLNL) hydrodynamic simulation code designed 
for large-eddy simulation of multicomponent flows with turbulent mixing~\cite{Coo99}. 
The spectral, high-order compact scheme with local artificial viscosity, 
diffusivity is used in order to remove oscillations and capture 
shocks and contact discontinuities~\cite{Coo07}. 
\textsl{RAGE} is a `radiation adaptive grid Eulerian' radiation-hydrodynamic code in which
the hydrodynamics is a basic Godunov solver~\cite{GitWeaClo09}.
\textsl{RAPTOR} is a hydrodynamic code 
based on a Godunov-type finite volume method that solves 
Riemann problem at cell interface using an adaptive mesh refinement technique~\cite{RanNieOak05}.
\textsl{TURMOIL} is a Lagrange-remap hydrocode which calculates the mixing of compressible fluids 
by solving the Euler equations plus advection equations for fluid mass fractions~\cite{You13}.
\FT\ is based on front tracking method for
accurate representation of the interface~\cite{GliGroLi98}. 
In this paper, we present the algorithmic features of \FT\ that are used to study the 
two-dimensional RMI instabilities. 
The front tracking algorithm is a way to track the interface explicitly with high order accuracy. 
It is a unique method demonstrated to avoid systematic errors in an important 
class of problems revolving around turbulent mixing \cite{GliCheSha20, GliShaKam11, ZhaKamShe18}.
This technique stores and dynamically evolves a meshed front that partitions 
a simulation domain into two or more regions, each representing a different 
material or physics model (see section~\ref{sec:ft}).

The incident shock strength characterized by Mach number has a big effect 
on the dynamics of flows.  
Some numerical methods such as filtered spectral methods only 
show numerically stable solutions under the moderate Mach number. 
When the Mach number is large, the methods become non-robust.  
Robust numerical methods based on solution-averaged or solution-reconstruction methods also known as reconstruction-evolutionary methods are 
used in the hydrodynamically unstable flow problems. 
Laney~\cite{Lan98} reviews highly used reconstruction-evolutionary methods including Van Leer's MUSCL scheme~\cite{Lee79},
Collella and Woodward's piecewise parabolic method (PPM)~\cite{ColWoo84}, 
Harten and Osher's uniformly high-order accurate non-oscillatory (UNO) schemes~\cite{HarOsh87}  
and
Harten, Engquist, Osher and Chakravarthy's essentially non-oscillatory (ENO) scheme~\cite{HarEngOsh87}.
ENO shock-capturing schemes were first introduced by Harten, Osher and 
Engquist~\cite{HarOshEng86}.  
The efficient implementation of ENO schemes based on numerical fluxes and total variation diminishing (TVD) Runge-Kutta 
time discretization are presented in \cite{ShuOsh89}. 
Liu, Osher and Chan \cite{LiuOshCha94} proposed the WENO schemes to shock capturing using a convex combination 
of all possible stencil candidates. 
Their WENO schemes use logical statements to choose ENO stencils which make the computation expensive. 
Jiang and Shu \cite{JiaShu96} introduced a new smoothness indicator and increased the fourth order scheme of Liu to the fifth order. 
Balsara and Shu \cite{BalShu00} extended it to higher-order.  
In the last two decades,  the weighted essentially non-oscillatory (WENO) methods gained popularity for solving hyperbolic systems of conservation laws to get accurate and robust solutions at sharp gradient regions~\cite{Shu20}. 
We use the classical fifth order WENO scheme of Jiang and Shu~\cite{JiaShu96} 
with Yang's artificial compression~\cite{Yan88} on the shock computation.

The main purpose of this paper is to predict the growth rate of the impulsively 
accelerated interface before and after reshock for two induced shock strengths. 
The organization of the paper is as follows. 
In section~\ref{sec:numerics} the numerical approaches for turbulence 
simulation and modeling approaches, the front tracking method and 
the fifth order WENO scheme are briefly reviewed. 
In section~\ref{sec:vvRTI} we present previous successful validation and verification studies of RTI using front tracking. 
In section~\ref{sec:results} we present the numerical results of the front tracking method with 
the fifth order WENO scheme on test problems.
We first introduce the two-dimensional Euler equations for inviscid flows and the
mass fraction equation in section~\ref{sec:Euler}. 
The numerical validation of WENO scheme 
on Sod's shock tube and Shu-Osher's shock-entropy wave interaction problems are presented in sections~\ref{sec:Sod1d} and \ref{sec:ShuOsher1d}.
The numerical solutions obtained using the classical fifth order WENO 
scheme of Jiang and Shu \cite{JiaShu96} 
with the artificial compression idea of Yang~\cite{Yan88} 
are compared with the reference solutions. 
In section~\ref{sec:vvRMI} we present our validation and verification results for the 
two-dimensional shock-driven Richtmyer--Meshkov Instability of an air/SF$_6$ interface problem.
We compare RMI simulations with the shock tube experiments conducted
by Collins and Jacobs.
For validation of RMI simulations,  in section~\ref{sec:validationRMI}
the single-mode shock tube experiments of Collins and Jacobs
with two incident shock waves ($M=1.11$ and $M=1.21$)~\cite{ColJac02} are introduced.
For verification of RMI simulations,  in section~\ref{sec:verificationRMI} 
we present the effect of mesh resolution 
on the quantities of interest such as interface displacement and amplitude.
We present the concluding remarks and discussion on the future improvement 
of RMI simulations in section~\ref{sec:conc}.

\section{Numerical Models and Algorithms}
\label{sec:numerics} 

Numerical methods for simulating turbulent mixing due to RTI have been widely explored. 
Schilling's paper~\cite{Sch20} on understanding RTI between simulation, modeling and experiment 
reviews the efforts.
In this section, a brief summary of numerical approaches used for modeling turbulence,  the front tracking method
and the WENO scheme is presented. 

\subsection{Turbulence Simulation and Modeling Approaches}
\label{sec:turb}
The hydrodynamic unstable turbulent flows such as KHI, RTI and RMI are deeply multiscale. 
The level of temporal and spatial scales of a turbulent flow
that are desired to be resolved identify the characteristic properties of the numerical approach. 
The mesh resolution and time steps required to solve the fluid structure depends on the Reynolds number,
which indicates the laminar, transitional, or fully turbulent flow.  

{\it Direct numerical simulation (DNS)} approach solves the full NSE 
numerically without modeling turbulence. 
All the temporal and spatial scales of the flow are resolved 
from the smallest scale (Kolmogorov's scale $\eta$) 
to the largest scale (integral scale $L$).  
DNS requires a mesh resolution that makes the computation 
unfeasible for engineering applications because the computational 
cost increases as the cube of the Reynolds number. 
While it is the most accurate numerical approach, 
it is limited to the low and moderate Reynolds numbers of flows and simple geometries.

{\it Reynolds-Averaged Navier-Stokes (RANS)} approach solves the time-averaged equations for mean values of all quantities. 
The additional nonlinear stress terms in momentum equation are modeled for the time-averaged solutions to the NSE.
RANS resolves length scales sufficient to specify the problem geometry. 
RANS is known as the most computationally efficient and common approach in engineering applications 
but it is the least accurate approach.  

{\it Large-Eddy Simulation~(LES)} approach solves the filtered NSE and resolves 
certain length scales that contain the most energy. The smaller scales are modeled using
the explicit subgrid scale models~\cite{GerPioMoi91,MoiSquCab91}.  
LES was first proposed by Smagorinsky \cite{Sma63} for the study of the 
dynamics of the atmosphere's general circulation. 
The multi-species compressible Navier-Stokes equation, 
filtered at a grid level, is solved so that the LES defines SGS terms 
(such as the Reynolds stress) as a source. 
These source terms are modeled as gradient diffusion terms, and 
the coefficients (turbulent viscosity, {\it etc.}) 
are recovered in a dynamic manner from the solution itself. 
That is called dynamic SGS model, the coefficients 
are computed locally in the simulation \cite{GerPioMoi91}. 
A broad problem, across many application areas, is to assess 
the accuracy of the SGS models.  

{\it Implicit Large-Eddy Simulation~(ILES)} approach solves the non-dissipative filtered NSE using the implicit 
subgrid scale model.  It models the effect of the small and unresolved scales by adding dissipation 
in the high wave number range. This approach limits the dissipation with the choice of the subgrid scale coefficients to 
constant numerical values from discretization errors or the nonlinear flux limiters~\cite{GriMarLG07}. 

\subsection{Front Tracking}
\label{sec:ft}
This is an adaptive computational method that provides sharp
resolution of a wave front by tracking the interface between distinct materials.  
Front-tracking (FT) discretizes and advects a surface representing the material interface 
which can then be used to adjust solution steps, so that the unwanted mixing between neighboring 
cells of different materials is prevented.
It represents interface explicitly as lower dimensional meshes moving through a rectangular grid. 
In 2D, the wave is represented by a curve which is comprised of connected line segments. 
In 3D, the wave is represented by a triangular mesh. 
The states (density, pressure and velocity) of fluids
are located in the centers of each grid cell~\cite{GliGroLi98, GliGroLi00}. 
The method solves the equations with the following main steps:
1. interface propagation,
2. interpolation reconstruction, and 
3. interior states update.
In the interface propagation step, the speed for each 
interface point is solved by either a
local Riemann problem or an interpolation method. Then, a new position for the point is
determined from the equation $x^{new}=x^{old} + V \Delta t {\mathbf n}$,
where V is the wave speed,  ${\mathbf n}$ is the normal direction on the point and 
$\Delta t$ is the time step size.
After all points in the interface are propagated, we get a new interface at a new time step. Methods are implemented in \FT\ to resolve the topological change and to optimize all triangles in the interface. 
In the interpolation reconstruction step, the components (defined
in the sense of point set topology relative to the ambient space with the interface removed)
at the cell center of all interior points are determined from the propagated interface. 
In the interior state update step, a Strang splitting is applied and three 1D equations are solved
consecutively. All the states on cell centers are updated using WENO scheme. If the stencil of
the scheme does not cross any interface, the states in the stencil are given by the cell center
values. Otherwise, a ghost cell method~\cite{GliIsaMar81} is used to fill the states on the points on the
other side of the interface.
\FT\ code has been developed for over three decades, and has been under continuous enhancement by the computational and applied mathematics research groups at Stony Brook University and the University of Arkansas. 

\subsection{WENO scheme}
\label{sec:weno}
The stable and high order WENO (Weighted Essentially Non-Oscillatory) scheme 
is widely used for turbulent mixing
due to Richtmyer--Meshkov Instability~\cite{LatSch20, MorSch14, MorSch13, SchLat10, SchLatDon07, LatSchDon07a, LatSchDon07b}. 
The main features of the WENO schemes designed to solve the 
hyperbolic equations in computational fluid dynamics can be found in \cite{Shu20}.
In 2D, the fluxes in $x$ and $y$ are calculated separately. 
High order accurate and non-oscillatory scheme flux reconstruction uses a convex combination 
of $k$ candidate stencils~\cite{JiaShu96}. 
For $k=3$, the fifth ($2k-1$) order finite difference WENO schemes approximates the derivative $ F({U} )_x$ at a point $x_i$,  
\begin{equation} 
{F} ({ U} )_x \vert_{x=x_i} \approx \frac{1}{\Delta x} ( {\hat F}_{i+1/2} - {\hat F}_{i-1/2} ) 
\end{equation} 
where ${U}$ is the state vector, ${F}({U}) $ is the flux,  
${\hat F}_{i+1/2}$ and  ${\hat F}_{i-1/2}$ are the fluxes at the cell boundaries,  and 
${\hat F}^+_{i+1/2}, {\hat F}^-_{i+1/2}$ are the positive and negative parts of  
$${\hat F}_{i+1/2} = {\hat F}^+_{i+1/2} + {\hat F}^-_{i+1/2}.$$
The fifth order WENO scheme uses three stencils 
$${\hat F}_{i+1/2}=\sum_{j=1}^3 \omega_i {\hat F}^{(j)}_{i+1/2}$$
with three third order fluxes ${\hat F}^{(j)}_{i+1/2}$ and the nonlinear weights $\omega_i$. 

The flux-averaged WENO scheme uses Lax-Friedrichs flux splitting method.
For hyperbolic conservation equations,  the nonlinear part of WENO is carried out in (local) characteristic fields. 
The flux ${\hat F}_{i+1/2}$
uses the average state ${\bar U}_{i+1/2}$ at cell boundary $(i+1/2, j)$. 
The implementation starts with computing the average state ${\bar U}_{i+1/2, j}$, 
then the left and right eigenvectors and eigenvalues of the Jacobian $ F^\prime({\bar U}_{i+1/2, j})$ at the average state. 
We transform the conservative fields and physical fluxes onto the characteristic fields 
using the left eigenvector matrix. 
We perform Lax-Friedrichs flux splitting and use the WENO reconstruction procedure 
to compute the positive and negative fluxes in the characteristic fields. 
We transform back the characteristic fluxes to the physical fluxes using the right eigenvector matrix.
We perform the same steps for the other direction in $y$ using the average states ${\bar U}_{i, j+1/2}$ for two-dimensional problems.

\section{Validation and Verification}
\label{sec:vvRTI}
FT for accurate representation of an interface is essential to achieve successful validation and verification studies. 
Modeling hydrodynamics instabilities without FT shows nonphysical mixtures produced by the numerical diffusion if the mesh is not fine enough.
Using FT/LES/subgrid scale models ~\cite{GliCheSha20, ZhaKamShe18, GliShaKam11} showed excellent agreement between experiments and numerical simulations of RTI. 
In the context of the RTI, ~\cite{GliShaKam11} examined the 
much-debated questions of models for initial conditions 
for a set of Smeeton-Youngs experiments~\cite{SmeYou87}
and the possible influence of unrecorded long wavelength 
contributions to the instability growth rate.  
The spectral amplitudes $A(n,t)$ for Fourier mode $n$ and
for time $t$ are used to estimate the spectral power law.  
A self-similar initial spectral power law 
$A(n)^2 \approx n^{-4}$ for the $t=0$ amplitudes was 
proposed \cite{You03,Dim04}. 
A similar power law was reported in the initial conditions 
for the surface of a glass Inertial Confinement Fusion (ICF) pellet, 
with higher amplitudes for the smaller wave lengths \cite{Dim04}. 
A power law $A(n)^2 \approx n^{-3.3}$ 
(mean exponent averaged over five experiments~\cite{SmeYou87, BurSmeYou84}) 
in the reconstructed initial time $t=0$ spectra was found~\cite{KamGliSha10, KauKamYu11, KamKauGli11}. 
The error estimates in the reconstructed initial data $t=t_0$ were 
used to establish an uncertainty quantification interval for numerical 
simulations of the RTI growth rate parameter, $\alpha$. 
The studies showed that the long wave length initial perturbations can affect 
the growth rate parameter $\alpha$s with at most a 5\% effect relative to the experimental data.
In order to find the optimal growth rate for RTI simulations, 
the uncertainty quantification studies based on 
the polynomial chaos expansion~\cite{Kam19} showed 
that the initial perturbation wavelength 
has an effect of 4\% on the quantity of interest,  
the growth rate parameter $\alpha$. 
The global sensitivity analysis 
addressed in \cite{Kam19} is consistent with \cite{GliShaKam11} that showed 
that the long wavelength initial perturbations can affect the growth rate parameter at most $\pm 5$ \%. 

The main focus of this paper is to address the issues in achieving agreement between simulations and experiment 
through validation and verification studies of  
single-mode RMI which are the shock-induced RTI. 
The robust front tracking method with WENO scheme is used for the 
numerical simulations of RMI of an air/SF$_6$ interface.
We investigate the dynamics of RMI with the single-mode perturbed interface 
and compare the numerical simulation results with the 
Collins and Jacobs~\cite{ColJac02} shock-tube experiments.
In the impulsively accelerated fluid-interface problems, 
the goal is to predict the amplitude and the displacement of the interface 
at all times including the pre-reshock and post-reshock.  
\section{Numerical Results}
\label{sec:results}

In section~\ref{sec:Euler} the two-dimensional (2D) Euler equations with the mass fraction solved for the inviscid numerical simulations are presented.
In sections~\ref{sec:Sod1d} and \ref{sec:ShuOsher1d} two benchmark test problems
i) Sod's shock tube and ii) Shu-Osher's shock--entropy wave interactions
are introduced for the numerical validation of the front tracking method with WENO scheme. 
In section~\ref{sec:vvRMI} we present our validation and verification results
for the single-mode shocked-induced RMI simulations of an air/SF$\_6$ interface.

\subsection{Euler equations}
\label{sec:Euler}
In 2D, the governing equations of the compressible inviscid gases are the Euler 
equations
\begin{equation}
{\bf U}_t + {\bf F} ( {\bf U})_x + {\bf G} ( {\bf U})_y={\bf 0}
\end{equation}
where ${\bf U}, {\bf F} ( {\bf U})$ and ${\bf G} ( {\bf U})$ are the vectors of 
conserved (mass, momentum, energy) variables and the fluxes in $x$ and $y$ direction.

\begin{equation}
\label{eq:Euler}
{\bf U}=\begin{pmatrix}
{\rho} \\
{\rho} {u} \\ 
{\rho} {v} \\ 
{E} \\
\end{pmatrix}, 
{\bf F} ( {\bf U})=
\begin{pmatrix}
{\rho} {u}\\
{\rho} {u}^2 + {p} \\
{\rho} u v \\
({E}+{p})u \\
\end{pmatrix}, 
{\bf G} ( {\bf U})=
\begin{pmatrix}
{\rho} v\\
{\rho} u v \\
{\rho} v^2 + p \\
({E}+{p})v \\
\end{pmatrix}.
\end{equation} 
Here $\rho$ is the density, $(u, v)$ is the velocity in $(x,y)$ directions, $p$ is the pressure, 
$E=\rho e + \frac{1}{2} \rho (u^2+v^2)$ is the total energy
$e=\frac{p} {(\gamma - 1) \rho }$ is the specific internal energy
$\gamma$ is the constant specific heat ratio.

The Euler equations~(\ref{eq:Euler}) are extended with the mass fraction equation
for cases with binary fluid mixing
\begin{equation}\label{eq:Conc}
\frac{\partial \rho M}{\partial t} + 
\frac{\partial \rho M u}{\partial x} + 
\frac{\partial \rho M v}{\partial y}
=0,
\end{equation}
where $M$ is the mass fraction for the heavy fluid.

\subsection{Sod's Shock Tube Problem}
\label{sec:Sod1d}
\begin{figure}[!ht]
    \begin{center}
       	 \includegraphics[width=.4\textwidth]{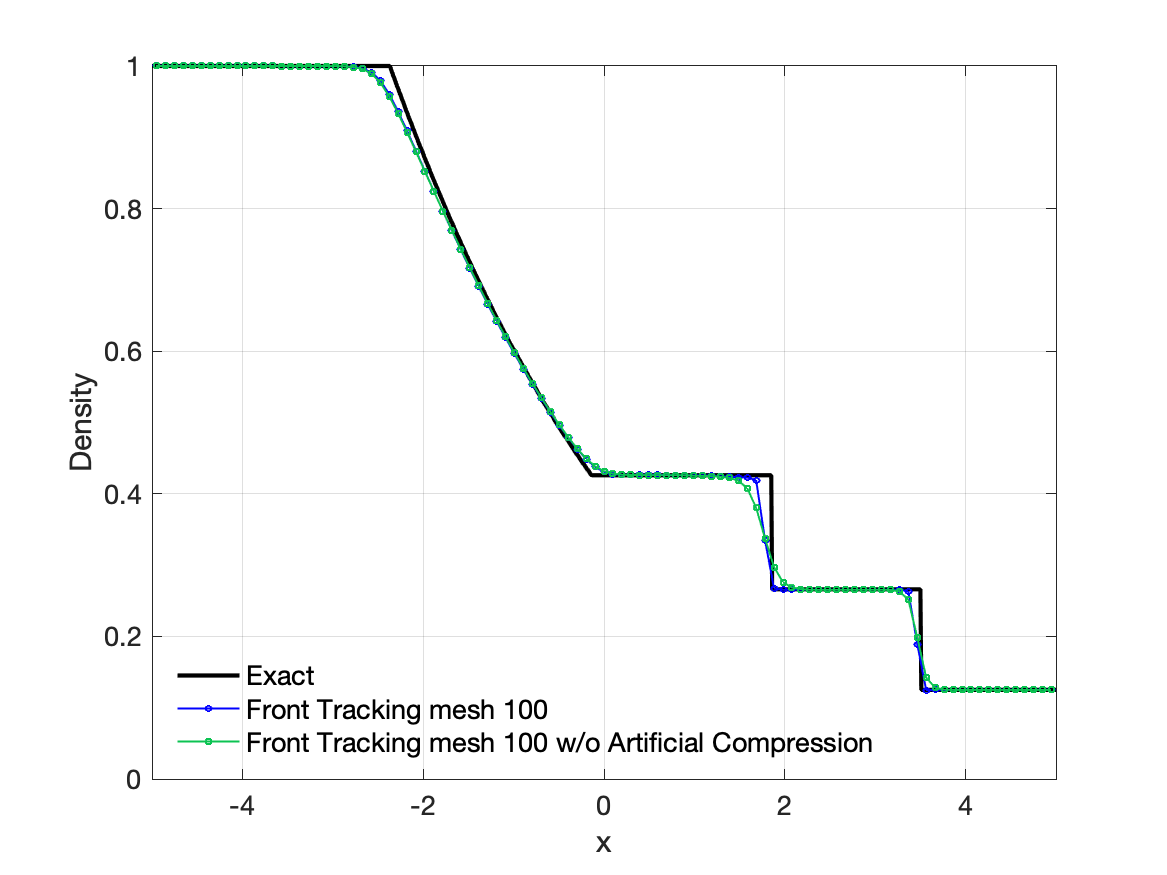}
     	 \includegraphics[width=.4\textwidth]{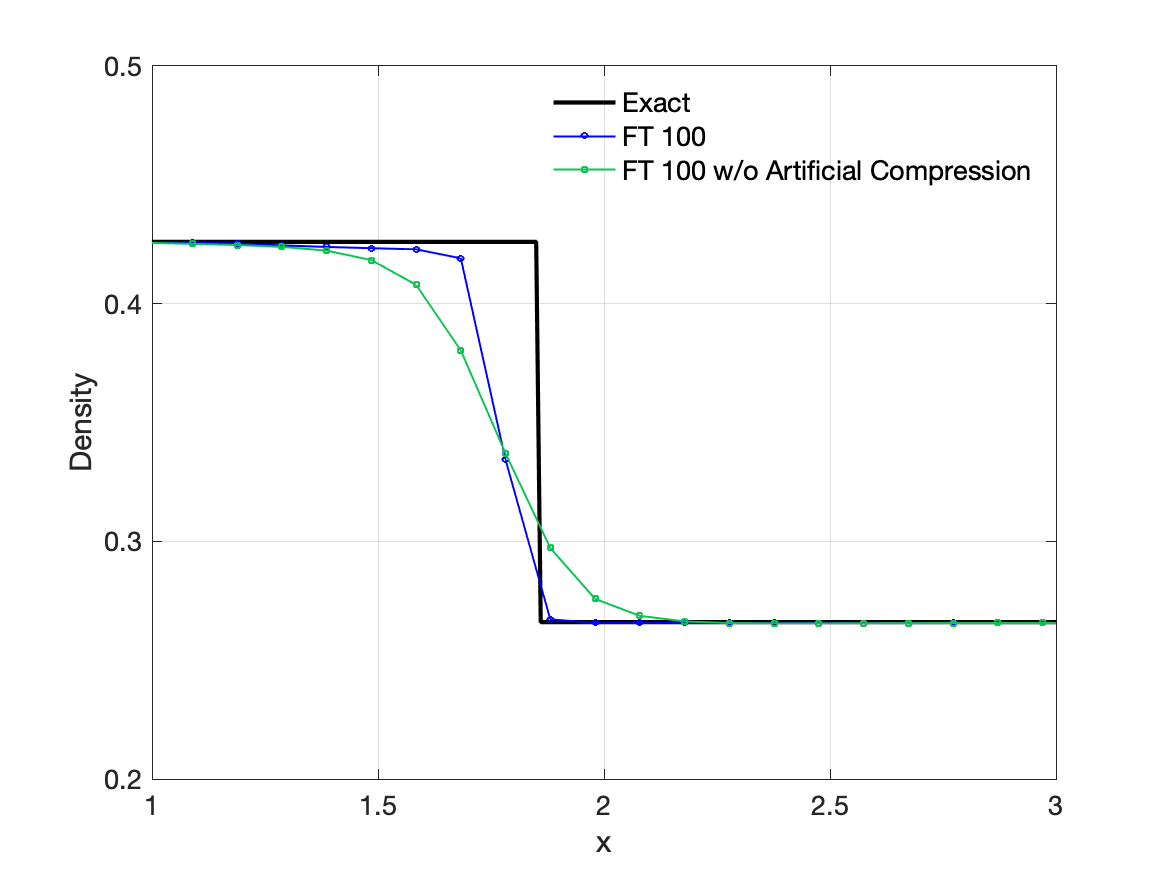}
 		 \caption{Left: The fifth order WENO scheme with and without artificial compression method of Yang~\cite{Yan88}
 		 are compared with the exact solution of Sod's shock tube problem. 
 		 Right: A zoomed view on domain $[1, 3]$. }
		\label{fig:1DsodShocknoA}
	\end{center}
\end{figure}

\begin{figure}
    \begin{center}
       	 \includegraphics[width=.4\textwidth]{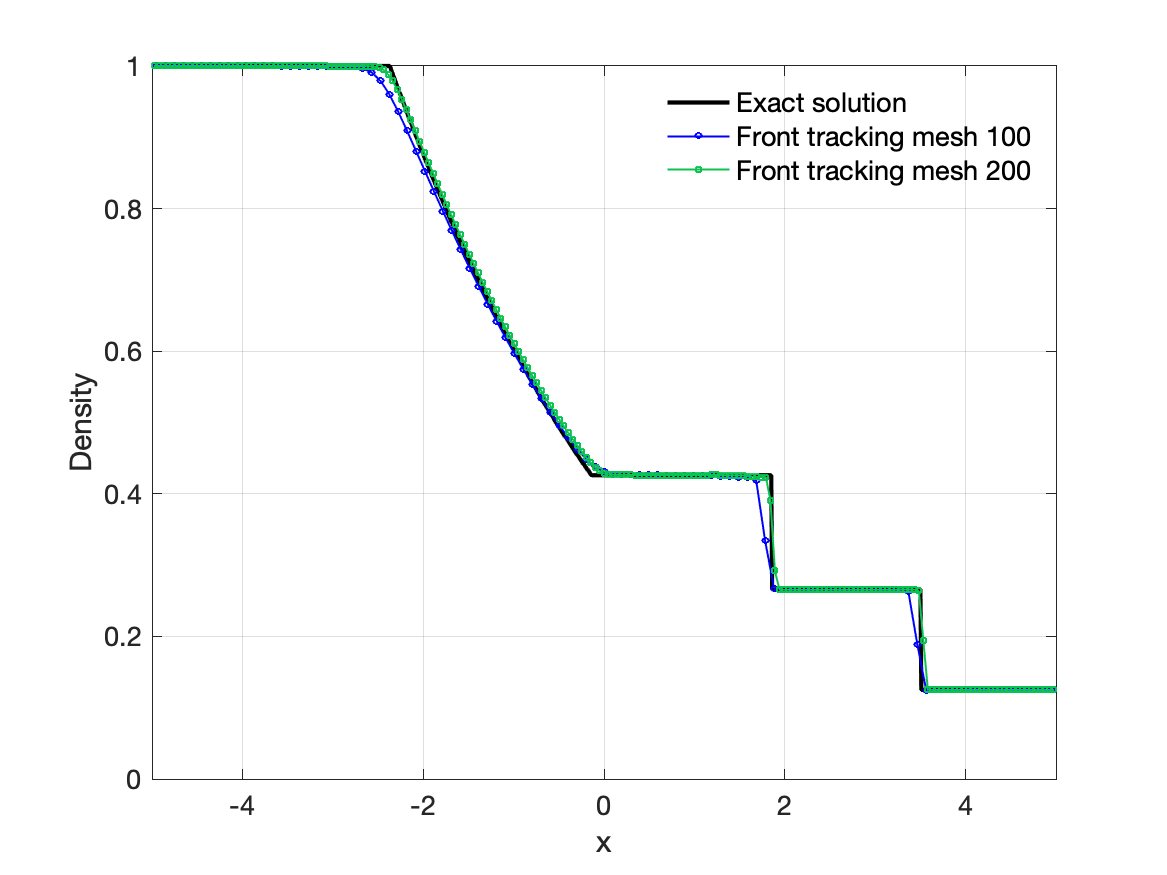}
     	 \includegraphics[width=.4\textwidth]{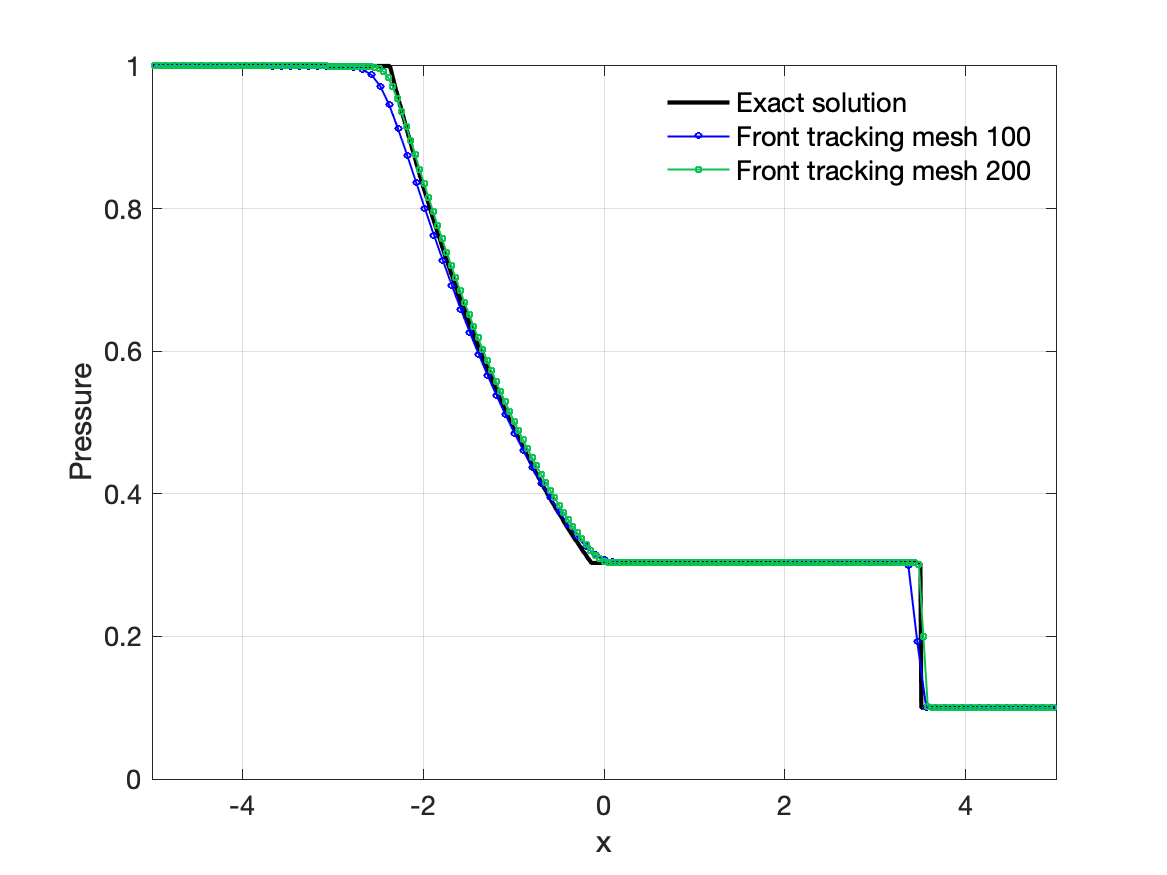}\\
     	  \includegraphics[width=.4\textwidth]{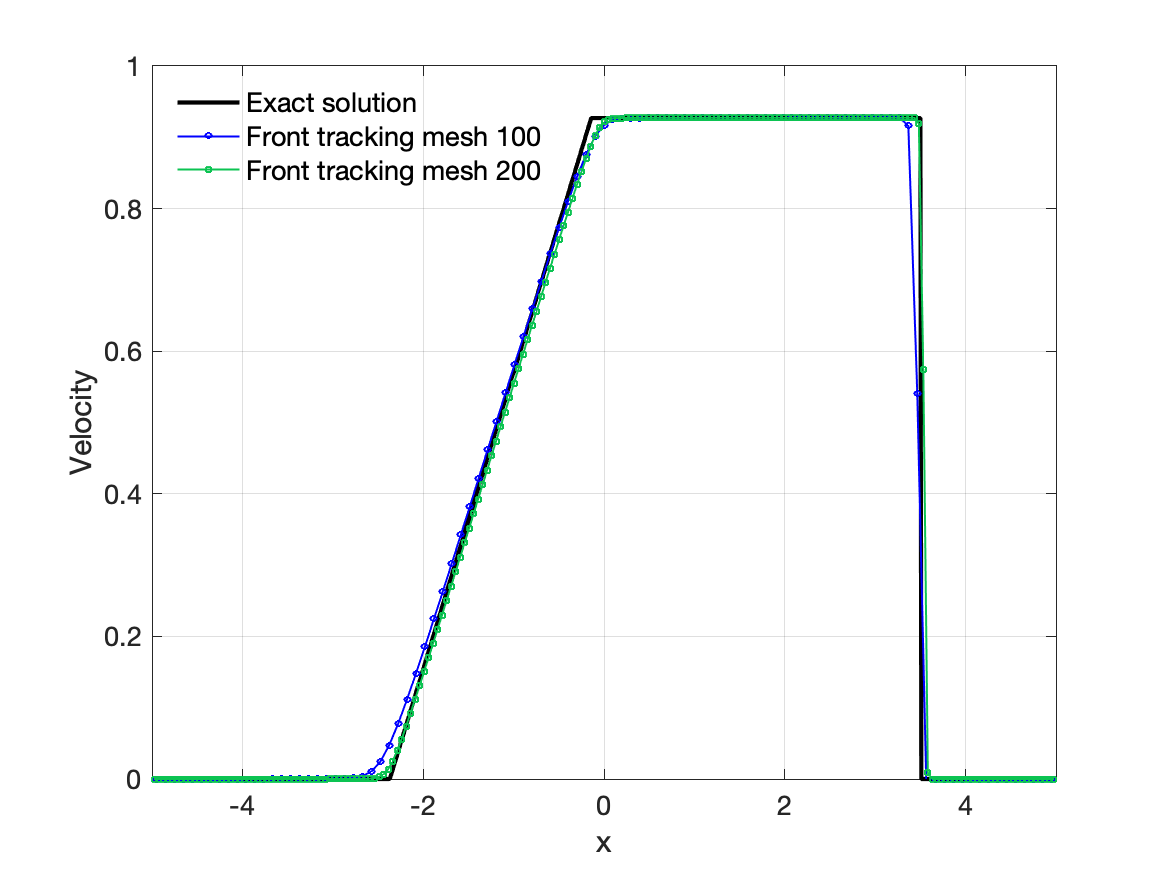}
 		 \caption{Density, pressure and velocity profiles on 100 and 200 mesh points are compared with the exact solution of Sod's shock tube problem.}
		\label{fig:1DsodShock}
	\end{center}
\end{figure}

For the contact discontinuity tracking,  
we consider the very well known shock-tube problem called Sod problem~\cite{Sod78} 
which is a Riemann problem with the initial condition as below.
\begin{equation}
  (\rho, u, p)=
    \begin{cases}
      (1, 0, 1) & \text{if $-5 \leq x \leq 0$}\\
      (0.125, 0, 0.1) & \text{if $0 \leq x \leq 5$}\\
    \end{cases}       
\end{equation}
The solution of Sod's shock tube problem contains
a rarefaction wave going to the left,
a contact discontinuity traveling slowly to the right, and 
a shock wave moving fast to the right.
The fifth order WENO scheme 
with and without Yang's artificial compression method on 100 mesh points is 
compared with the exact solution.  The slope correction in Yang's artificial compression method 
shows improvement at the contact discontinuity (see Figure~\ref{fig:1DsodShocknoA}).
In Figure~\ref{fig:1DsodShock} the density, pressure and velocity profiles 
with artificial compression on 100 and 200 mesh points
at $t=2.0$ are compared with the exact solution. 

\subsection{Shock--entropy wave interaction}
\label{sec:ShuOsher1d}
\begin{figure}
\begin{center}
  \includegraphics[width=.42\textwidth]{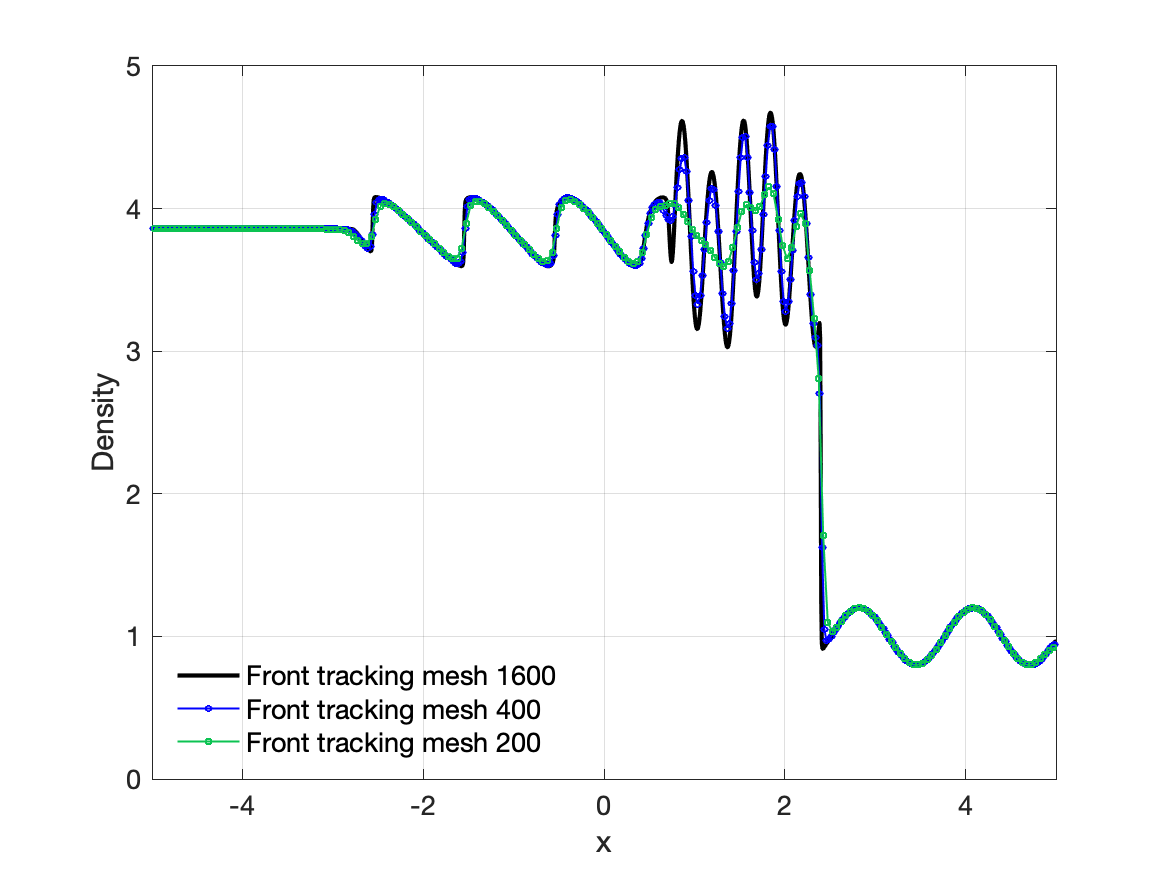}
  \includegraphics[width=.42\textwidth]{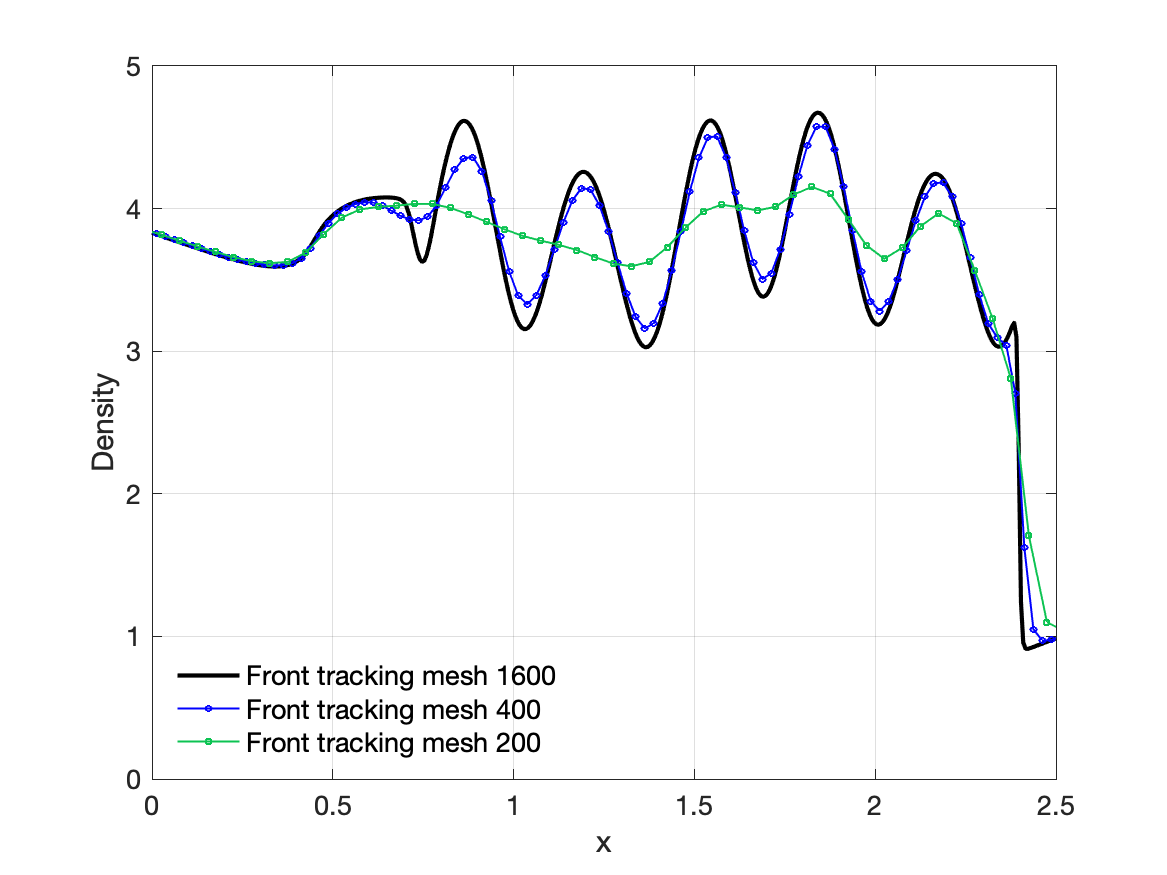}
  \caption{Left: Density profile of the shock-wave interaction problem 
  on domain $[-5,5]$ 
  with $200,400$ and $1600$ mesh points.
  Right: A zoomed view on domain $[0,2.5]$. }
  \label{fig:1Dshockentropy}
  \end{center}
\end{figure}

In order to study the stability and accuracy of the WENO scheme for 
the interaction between a shock wave and an entropy wave, 
we consider the Riemann problems for the Euler equations 
for an ideal gas with a constant ratio of specific heats 
($\gamma=1.4$). Shu-Osher's test problem~\cite{ShuOsh89} 
corresponds to a Mach $M=3$ shock wave passing 
through an entropy wave 
on the spatial domain $(-5, 5)$ and the time domain $(0,2)$. 
The shock wave and its interaction with sine waves in the density
field on 200, 400 and 1600 mesh points at $t=1.8$ are presented in 
Figure~\ref{fig:1Dshockentropy}.
\begin{equation}
  (\rho, u, p)=
    \begin{cases}
      (3.857143, 2.629369, 10.33333) & \text{if $x \leq -4$}\\
      ( 1 + \epsilon \sin kx, 0, 1) & \text{if $x \geq -4$}\\
    \end{cases}       
\end{equation}
where $\epsilon=0.2$ and $k=5$ are the amplitude and the wave number of the entropy wave.
The solution near N-wave ($x\in[-2,0]$), 
the transition to the N-wave ($x\in[0,1]$), 
the entropy wave ($x\in[1,2]$), and 
the shock ($x\in[2,2.5]$) is displayed on three levels of mesh refinement 
in Figure~\ref{fig:1Dshockentropy}.

\subsection{RMI simulations}
\label{sec:vvRMI}

In the impulsively accelerated RTI which is known as RMI,  the amplitude and displacement measurements are 
the quantities of interest.  We are interested in predicting the growth rate of the instability ($da/dt$) where $a$ is the amplitude 
and $t$ is the time.  Richtmyer~\cite{Ric60} modeled and defined the growth rate in terms of 
the wavenumber ($k=2\pi /\lambda$), where $\lambda$ is the perturbation wavelength,
the initial perturbation amplitude ($a_0^+$), 
the Atwood number ($A^+$), and 
the interface velocity ($\Delta u$) as below. 
\begin{equation}
\frac{da}{dt} = V_{RM} = k a_0^+ A^+ \Delta u 
\label{eq:Ric60_growthRate}
\end{equation}
The Atwood number,  a dimensionless measure of the density contrast,  is defined by
\begin{equation}
A^+=\frac{\rho_2^+ - \rho_1^+}{\rho_1^+ +\rho_2^+}
\label{eq:Atwood}
\end{equation}
where $\rho_1^+$ and $\rho_2^+$ are the densities of the light and heavy fluids.  
The (+) sign indicates the post-shock quantity. 
\subsubsection{Validation}
\label{sec:validationRMI}
\begin{figure}
\begin{center}
  \includegraphics[width=.275\textwidth]{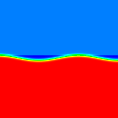}
  \includegraphics[width=.275\textwidth]{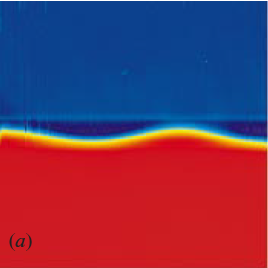}
  \caption{Initial diffusive sinusoidal perturbed interface in simulation (left) and experiment (right) before the shock wave.
  The light fluid (air) at the top and heavy fluid (SF$_6$) at the bottom.  
  Experimental figure courtesy of J.W. Jacobs. 
  The experimental image is taken from 
  Figure 6 of Collins and Jacobs (reprinted with the permission of Cambridge University Press).}
  \label{fig:initialPert}
  \end{center}
\end{figure}

For validation study,  the model and input parameters used in our RMI simulations 
are set according to Collins and Jacobs~\cite{ColJac02} shock tube experiments
conducted for the two shock wave Mach numbers $M=1.11$ and  $M=1.21$.  
We choose these experiments to compare with our numerical results because of the 
high quality images obtained using planar laser-induced fluorescence (PLIF)
visualization technique in the experiments.  The shock tube is 430~cm long 
and consists of a driver and a driven section which are 100~cm and 330~cm long 
respectively.  
The sinusoidal initial perturbation is created within the shock tube 
{\it test} section that has a 8.9 cm square cross section and a length 75~cm.  
In our numerical studies, the computations domain is set to 
$L_x=8.9$~cm and $L_y=75$~cm. 
The location of the initial shock ($y=1$~cm) 
and the perturbed interface ($y=3$~cm)
are set according to Latini, Schilling and Don 2007 numerical studies~\cite{LatSchDon07a, LatSchDon07b}.
The distance from the initial perturbed interface to the end of the shock tube 
is $72$~cm.
In Collins and Jacobs experiments 
the sinusoidal perturbed interface 
between the light fluid (air) and the heavy fluid (SF$_6$)
is initialized with an amplitude of $a_0=0.229$~cm for $M=1.11$ 
and $a_0=0.183$~cm for $M=1.21$. 
The diffuse interface thickness is set to $\delta=0.5$~cm  between 
the two fluids. 
The initial perturbation wavelength is $\lambda_0=5.933$~cm
across the width of shock tube.
In our numerical simulations, 
the initial density profile over the diffuse interface is set by
$\rho(y) = {\bar \rho} \left( 1 + A\  {\rm erf} (\sqrt{\pi} y / \delta) \right)$ 
where ${\bar \rho} = (\rho_1+\rho_2)/2$ is the average density, 
$A$ is the Atwood number and $\delta$ is the initial diffuse interface thickness.
Figure~\ref{fig:initialPert} shows the diffuse sinusoidal perturbed initial interface 
for the $M=1.21$ simulation at time $t=0.05$~ms and
the first experimental PLIF image taken just before the incident 
shock arrives to the interface. 
In this section, the fine grid simulations, that are compared with experiments,
have $256$ grid points per wavelength with a grid spacing
$\Delta x = \Delta y = 0.0232$~cm.
Table~\ref{table:params} shows the input parameter values including 
the density, molecular weight and the ratio of specific heats 
corresponding to the light and heavy fluids,  the pressure at the 
interface and the CFL number used in our numerical simulations.

There are some differences in the input parameters between our 
and Latini, Schilling and Don's numerical simulations 
(see~\cite{LatSchDon07a, LatSchDon07b} 
for more details on their specification of the input parameters).  
In \cite{LatSchDon07a, LatSchDon07b}, 
the length of the domain $L_y=78$~cm and the initial amplitude $a_0=0.2$~cm 
are larger than the originally reported experimental values. 
The simulations presented in here are set using exactly the same 
experimental parameter values reported in \cite{ColJac02} 
for the validation of our simulations. 
In terms of repeatability of experiments, Collins and Jacobs reported 
10\% difference in the initial perturbations. 
The effect of the input parameters on the quantities of interest 
such as amplitude and displacement will be addressed in the future publication.

\begin{table}
\begin{tabular}{ p{7cm}p{2cm}p{1cm}}
\hline
 \multirow{2}{*}{} &  \multicolumn{2}{c}{Mach Numbers}  \\ 
                              & $1.11$ &  $1.21$ \\ \hline
 Initial amplitude $a_0$~(cm)   & 0.229 & 0.183 \\
 Initial wavelength $\lambda_0$~(cm)  & \multicolumn{2}{c}{5.933} \\
 Ratio between $a_0/\lambda_0$  & $0.0386$ & $0.0308$  \\
 Initial diffusion layer $\delta$~(cm)  & \multicolumn{2}{c}{0.5} \\
 Heavy fluid (SF$_6$) density $\rho_1$~(g/cm$^3$)   & \multicolumn{2}{c}{$5.494 \times 10^{-3}$} \\
 Light fluid (air) density $\rho_2$~(g/cm$^3$)  & \multicolumn{2}{c}{$1.351 \times 10^{-3}$} \\
 Atwood number $A$   & \multicolumn{2}{c}{0.6053} \\
 Molecular weight of SF$_6$~(g/mol)  & \multicolumn{2}{c}{$146.05$} \\
 Molecular weight of air~(g/mol)   & \multicolumn{2}{c}{$34.76$} \\
 Ratio of specific heats $\gamma $   & \multicolumn{2}{c}{1.276} \\
 Pressure at interface $p$~(bar)  & \multicolumn{2}{c}{0.956} \\ 
 Courant--Friedrichs--Lewy (CFL) number   & \multicolumn{2}{c}{0.45} \\
 \hline 
\end{tabular}
\caption{The parameter values of RMI simulations.}
\label{table:params}
\end{table}

\begin{figure}
\begin{center}
  \includegraphics[width=.7\textwidth]{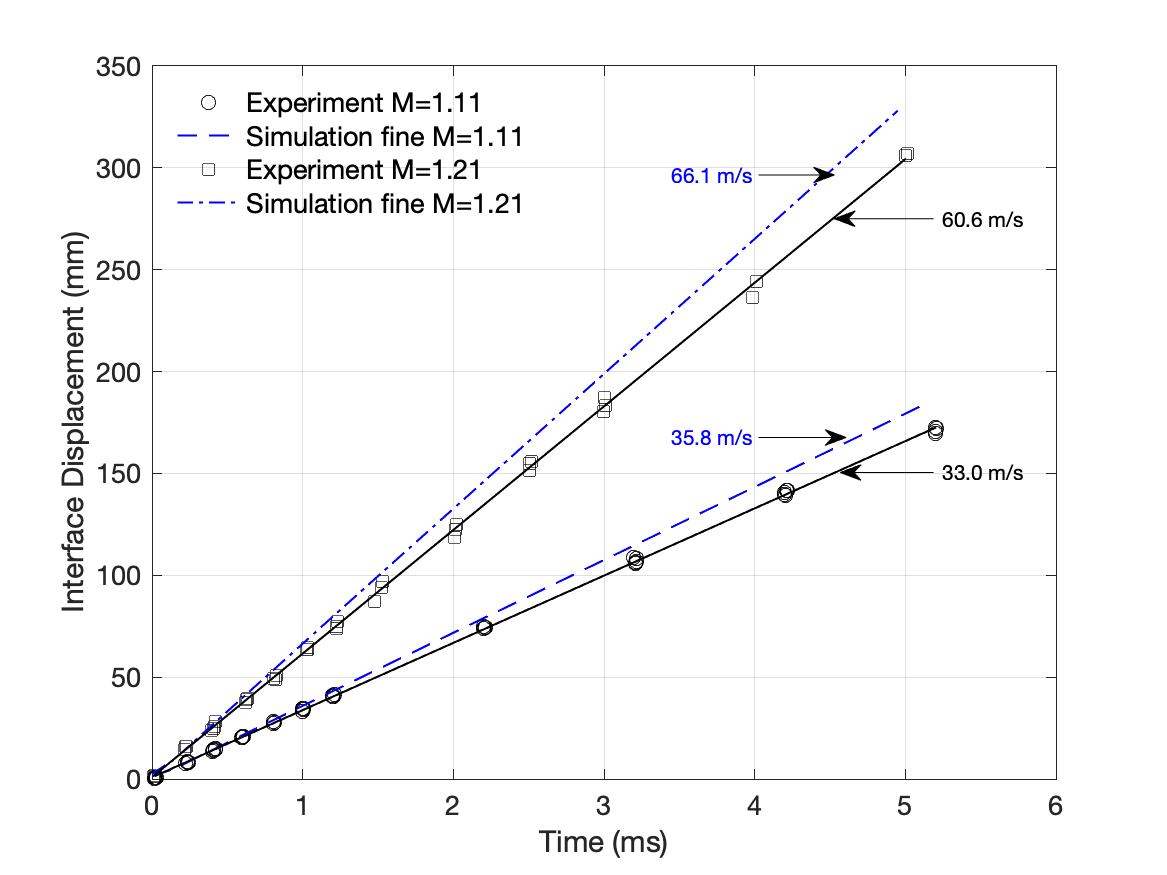}\\
  \caption{Comparison of interface displacement (mm) and interface velocity between Collins and Jacobs 2002 experiments and the
  \FT\ fine grid simulations for $M=1.11$ and $M=1.21$.}
  \label{fig:RMI2d_disp}
  \end{center}
\end{figure}

The light fluid (air) penetrates the heavy fluid (SF$_6$) in bubbles 
and the heavy fluid penetrates the light fluid in spikes.  
The interface displacement is computed by averaging the distance of the bubble and spike tip locations to the initial tip locations. 
In figure~\ref{fig:RMI2d_disp} the displacement of the interface is compared with the experiments for 
shock strengths $M=1.11$ and $M=1.21$. 
The values for the interface velocities $\Delta u$ of the simulations and experiments are given in Table~\ref{table:comp_exp_sim}.
The experimental measurements show a constant velocity of
$33.0$~m/s for $M=1.11$ and $60.6$~m/s for $M=1.21$, whereas 
simulation interface velocity is $35.8$~m/s for $M=1.11$ and $66.1$~m/s for $M=1.21$. 
Theoretical values computed using 
a one-dimensional analysis assuming ideal gas behaviour 
and a nonventilated shock tube with the measured shock velocities
are $36.0$~m/s for $M=1.11$ and $64.2$~m/s for $M=1.21$~\cite{ColJac02}. 
The values of the interface velocities of the numerical simulations 
are closer to the theory than the experiments.
It is reasonable to expect that the discrepancies in the interface velocities between the simulations and experiments 
are due to the two open slots in the test section of the shock tube experiments. 
When the incident shock passes through the slots in the experiments, it creates expansion waves
which slow down the interface velocity.  Because of the expansion waves which are only present 
in the experiments,  a difference less than 10\% is observed (see table
~\ref{table:comp_exp_sim}).

\begin{table}
\begin{tabular}{lp{1cm}p{1cm}p{1cm}p{1cm}p{1cm}p{1cm}p{1cm}p{1cm}p{1cm}p{1cm}p{1cm}p{1cm}}
\hline	
			       &  \multicolumn{6}{c}{M=1.11} &  \multicolumn{6}{c}{M=1.21} \\
\hline			       
                    &  \multicolumn{2}{c}{sim.} &  \multicolumn{2}{c}{exp.~\cite{ColJac02}} & \multicolumn{2}{c}{discrepancy} 
                    &  \multicolumn{2}{c}{sim.} &  \multicolumn{2}{c}{exp.~\cite{ColJac02}} & \multicolumn{2}{c}{discrepancy} \\
$\Delta u$(m/s) & \multicolumn{2}{c}{35.8} &  \multicolumn{2}{c}{33.0}  & \multicolumn{2}{c}{-8.5\%} 
						& \multicolumn{2}{c}{66.1} & \multicolumn{2}{c}{60.6}  & \multicolumn{2}{c}{-9.1\%} \\
$V_{RM}$(m/s)  & \multicolumn{2}{c}{4.10} &  \multicolumn{2}{c}{$3.92\pm0.23$}  & \multicolumn{2}{c}{0\%} 
& \multicolumn{2}{c}{5.44} & \multicolumn{2}{c}{$6.28\pm0.6$} & \multicolumn{2}{c}{4\%} \\
\hline
\end{tabular}
\caption{The simulation and experiment values for the interface velocity $\Delta u$ and Richtmyer velocity $V_{RM}$.
Discrepancy refers to the comparison of numerical values outside of uncertainty values, if any. 
}
\label{table:comp_exp_sim}
\end{table}

\begin{figure}[!ht]
\begin{center}
  \includegraphics[width=.475\textwidth]{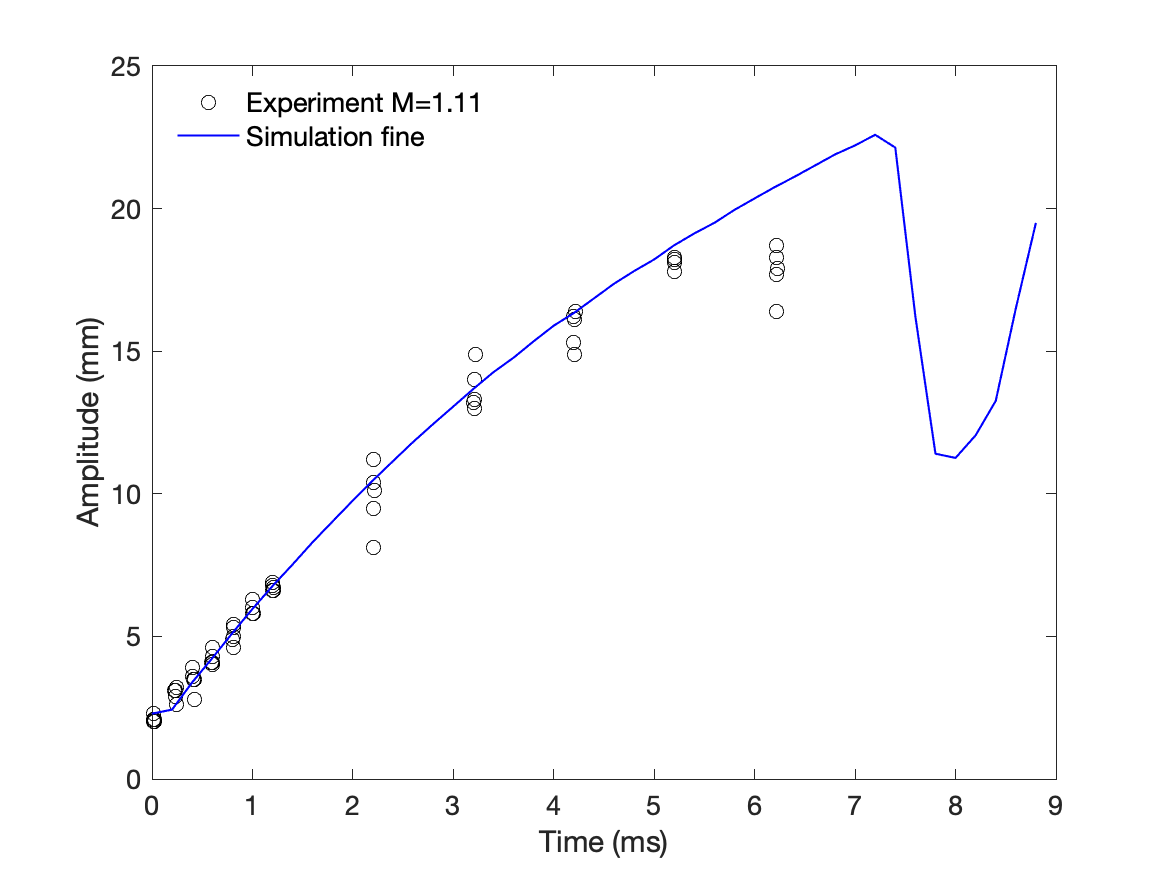}
  \includegraphics[width=.475\textwidth]{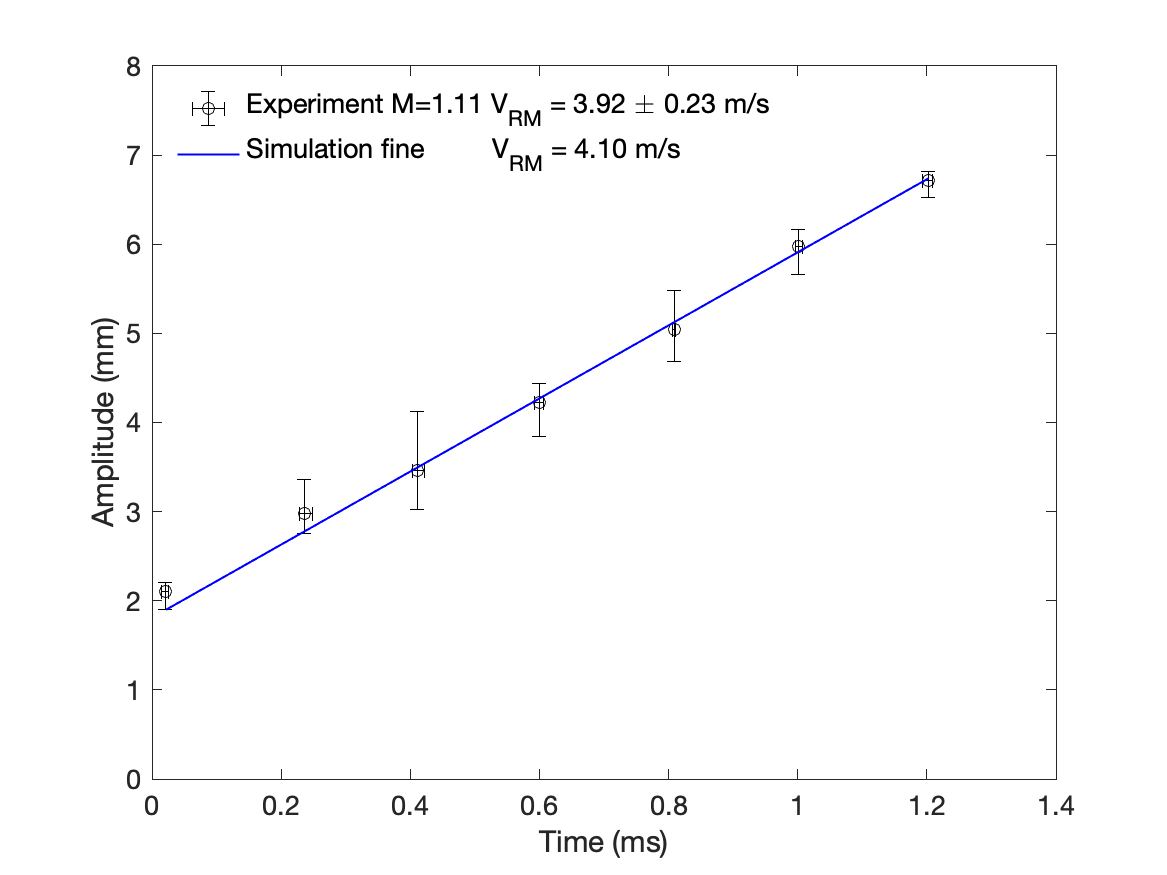}
  \caption{Comparison of amplitude between Collins and Jacobs 2002 
  experiments and \FT\ fine grid simulation for $M=1.11$. 
  Right: The early-time amplitude ($t \approx 1.2$ ms).
  Each experimental data point with error bar in time (ms) and amplitude (mm) shows the measurement from five experiments.}
  \label{fig:RMI2d_ampM1_11}
  \end{center}
\end{figure}

The amplitude is computed as half the vertical distance between the bubble and spike tip locations.
Figure~\ref{fig:RMI2d_ampM1_11} shows the amplitude data for the $M=1.11$ 
Collins and Jacobs 2002 experiments and the fine grid simulation.  The early-time amplitude ($t \approx 1$ ms) shows linear growth.
In figure~\ref{fig:RMI2d_ampM1_11} (right) each experimental data point 
with error bar in time (ms) and amplitude (mm) 
shows the measurement from five experiments.
The numerical simulation result $4.10$~m/s stays within the experiment 
error margin $3.92\pm0.23$~m/s for $M=1.11$.  
The shock is reflected at the end of the test section and the reflected reshock 
reaches to the interface at $t \approx 7.2$~ms, 
which is about $1$~ms later than the experiment. 

\begin{figure}[!ht]
\begin{center}
   \includegraphics[width=.49\textwidth]{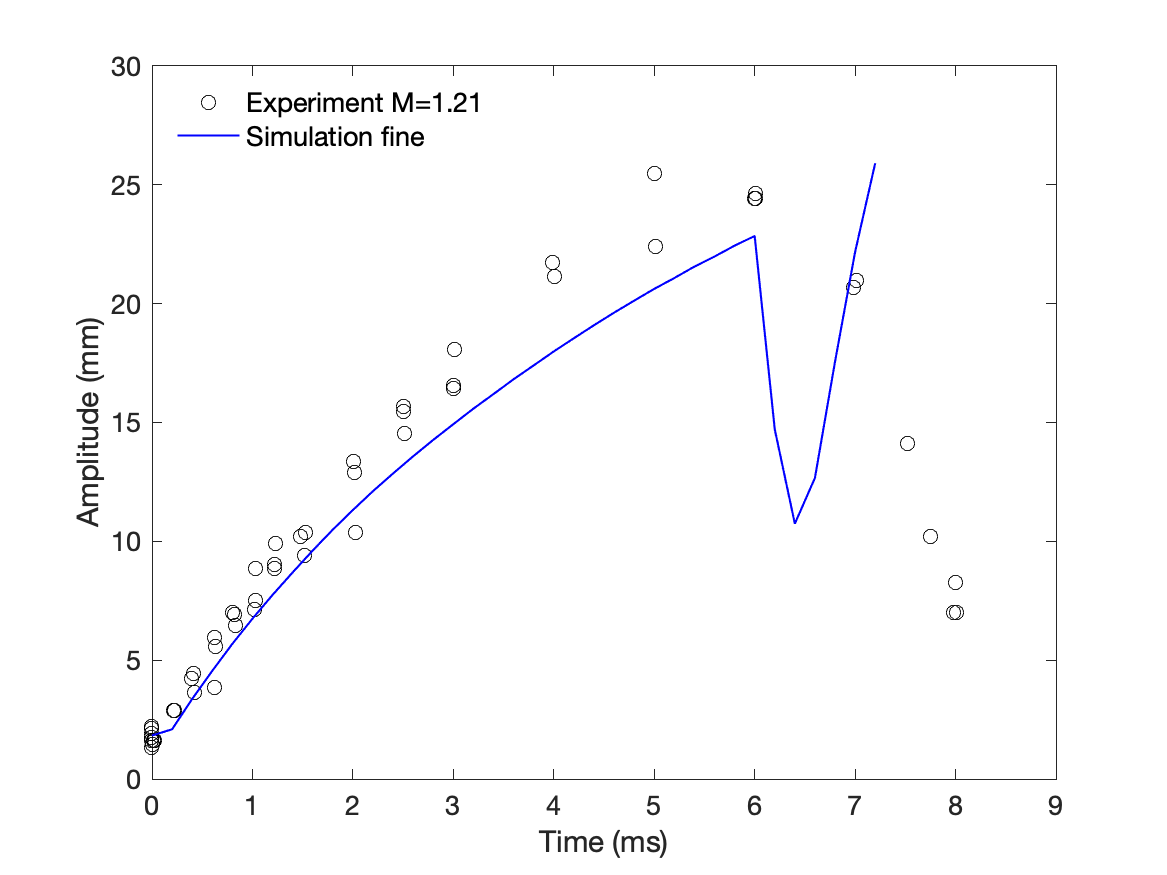}
  \includegraphics[width=.49\textwidth]{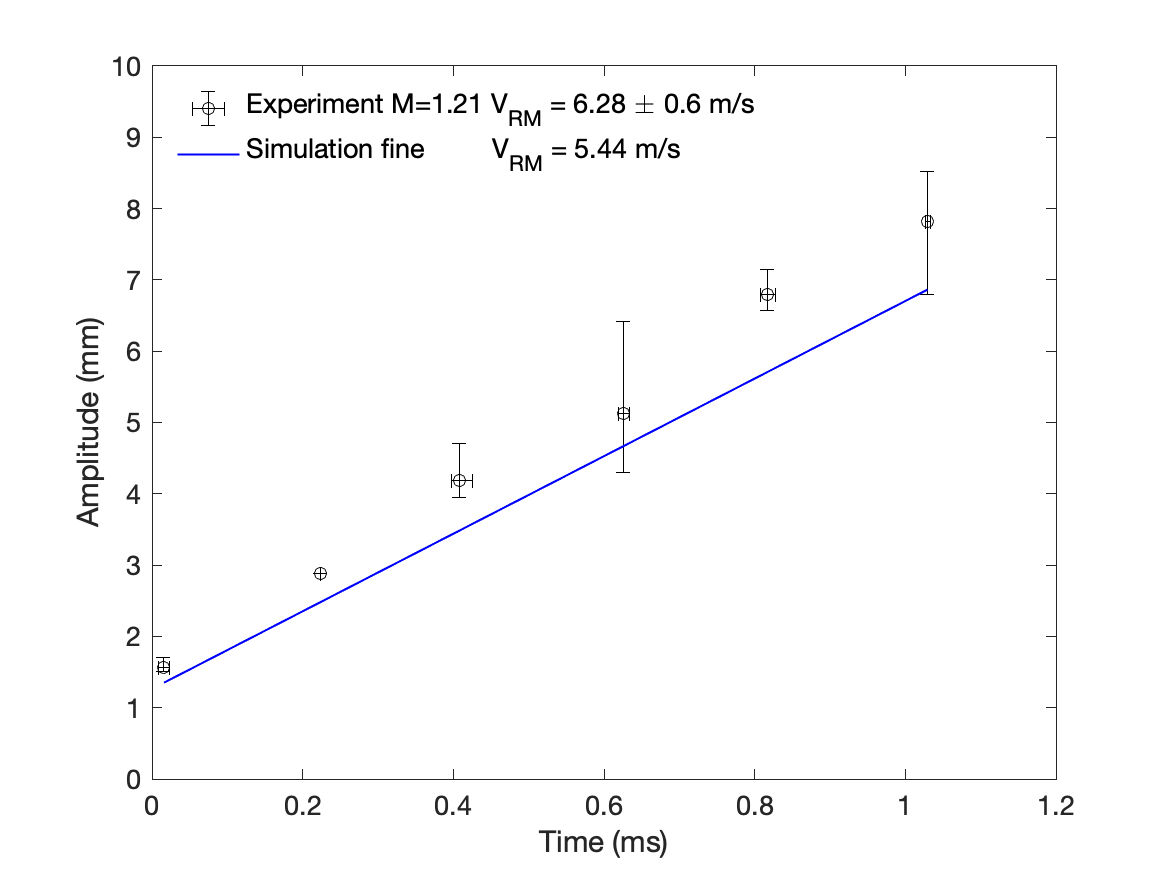}
  \caption{Comparison of amplitude between Collins and Jacobs 2002 experiments and \FT\ fine grid simulation for $M=1.21$. 
  Right: The early-time amplitude ($t \approx 1$ ms).
  Each experimental data point with error bar in time (ms) and amplitude (mm) 
  show the measurement from three experiments.}
  \label{fig:RMI2d_ampM1_21}
  \end{center}
\end{figure}

\begin{figure}
\begin{center}
 \includegraphics[width=.25\textwidth]{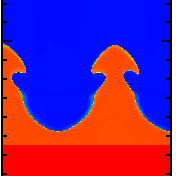}
     \includegraphics[width=.25\textwidth]{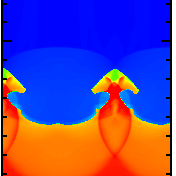}
           \includegraphics[width=.25\textwidth]{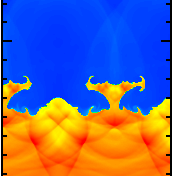}\\
   \includegraphics[width=.25\textwidth]{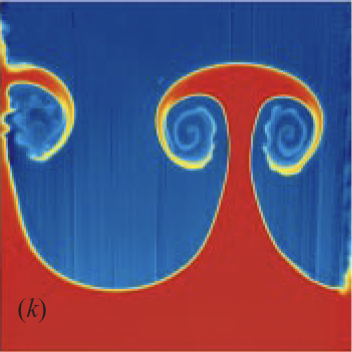}
     \includegraphics[width=.25\textwidth]{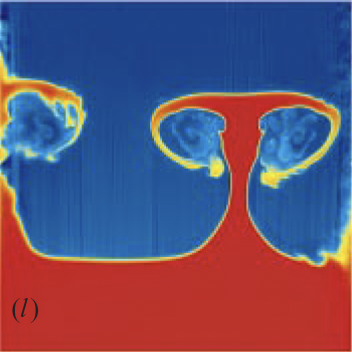}
       \includegraphics[width=.25\textwidth]{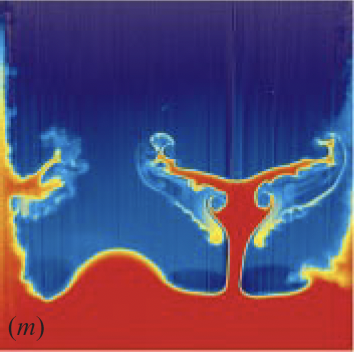}
  \caption{Top: The density profile of the fine grid simulation at times $6, 6.2$ and $6.4$~ms.  
  Bottom: Experimental PLIF images at times $6.006, 7.005$ and $7.781$~ms. 
  The experimental images are taken from 
  Figure 6 of Collins and Jacobs (reprinted with the permission of Cambridge University Press).}
  \label{fig:RMI2d_intfc_images}
  \end{center}
\end{figure}

Figure~\ref{fig:RMI2d_ampM1_21} shows the amplitude as a function of time for  
the $M=1.21$ Collins and Jacobs experiments and fine grid simulation. 
A 4\% discrepancy between the simulation and experiment early-time growth 
rate is observed.  In the present study,  the shock wave arrives to the perturbed interface at time $t=0.05$~ms (see figure~\ref{fig:initialPert}).  
The reflected shock reaches the interface at $t \approx 6$~ms (see figure~\ref{fig:RMI2d_intfc_images}). 
After reshock,  a time delay $\approx 1.2$~ms between simulations and experiments is observed 
because the expansion waves in the experiments reduce the interface velocity.
Therefore, in this present study 
the density profile from the experimental PLIF images at time $6.006$~ms
is compared with simulation data at time $6$~ms before reshock and 
after reshock PLIF images at times $7.005$~ms and $7.781$~ms 
are compared with the fine simulation density profile at times $6.2$~ms and $6.4$~ms. 
A similar time delay $\approx 1$~ms is also observed by Latini, Schilling and Don 2007~\cite{LatSchDon07b}.  
This time difference between numerical simulations is because of the difference in the length of the domain 
between our and Latini, Schilling and Don 2007 simulations. 
Figure~\ref{fig:RMI2d_intfc_images} shows a good agreement on the interface 
overall structure, but some of the finer details such the vortex roll-up structure
are missing. 
The missing issue in the transport of vorticity along the interface requires an 
upgrade in the \FT\ and WENO code.
\subsubsection{Verification}
\label{sec:verificationRMI}
\begin{figure}
\begin{center}
  \includegraphics[width=.75\textwidth]{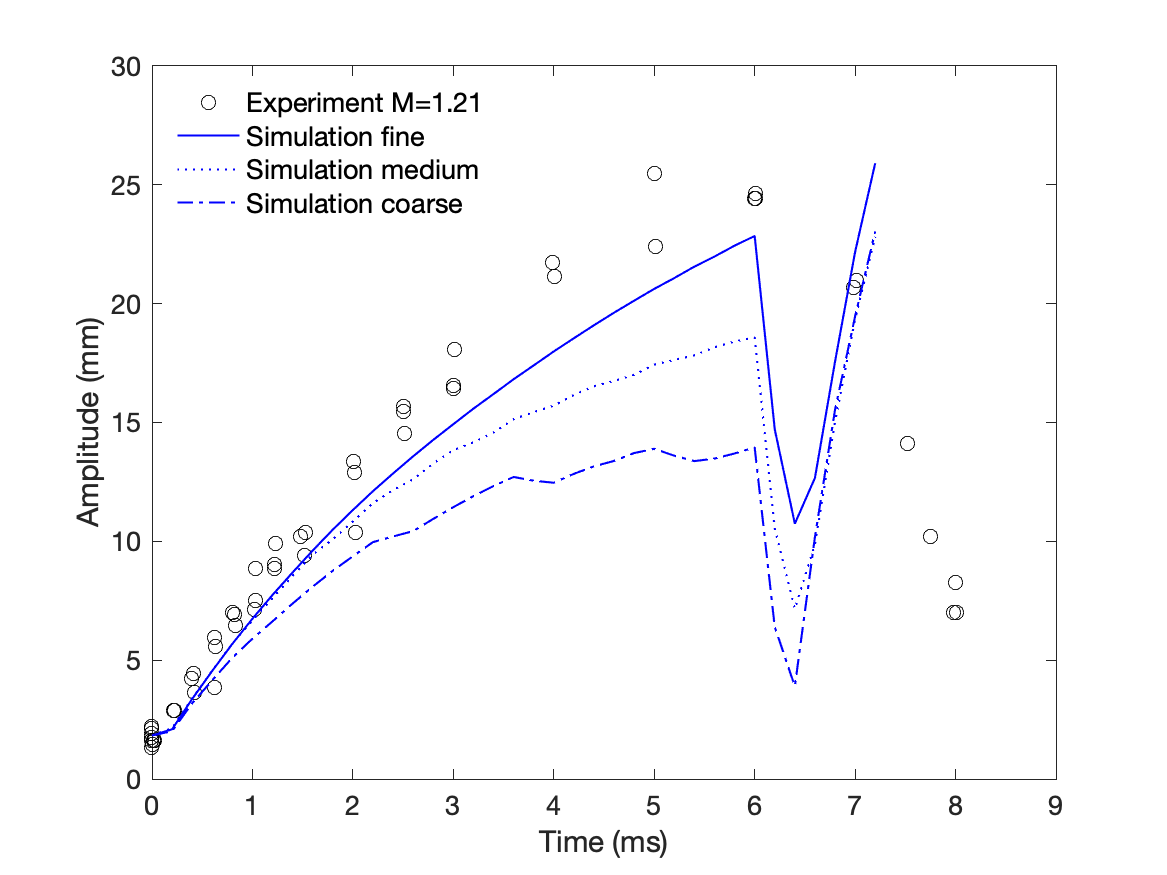}    
  \caption{Comparison of the amplitude between Collins and Jacobs 2002 experiments and \FT\ simulations for $M=1.21$.  
  The number of grid points per initial perturbation wavelength is 
  $64, 128$ and $256$ for the coarse, medium and fine grid simulations.}
  \label{fig:amp_sim3}
  \end{center}
\end{figure}

\begin{figure}
\begin{center}
  \includegraphics[width=.25\textwidth]{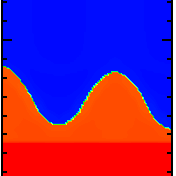}
     \includegraphics[width=.25\textwidth]{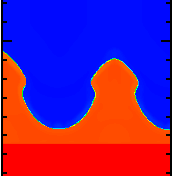}
      \includegraphics[width=.25\textwidth]{figs/II-f-FT-dens-M1_21-t6.png} \\
  \includegraphics[width=.25\textwidth]{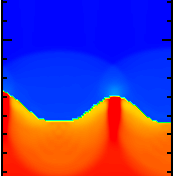}
     \includegraphics[width=.25\textwidth]{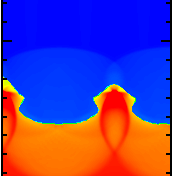}
      \includegraphics[width=.25\textwidth]{figs/II-f-FT-dens-M1_21-t6_2.png} \\
        \includegraphics[width=.25\textwidth]{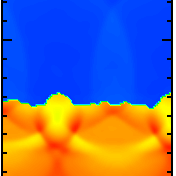}
     \includegraphics[width=.25\textwidth]{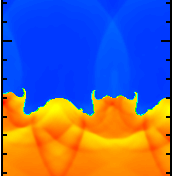}
      \includegraphics[width=.25\textwidth]{figs/II-f-FT-dens-M1_21-t6_4.png}
  \caption{Left to right: The density profile from coarse to fine grid simulations.
  Top to bottom: At times 6, 6.2 and 6.4~ms. }
  \label{fig:fmc_densProfile}
  \end{center}
\end{figure}

While validation of the numerical simulations is done through comparison with experiments,
verification is done through mesh refinement. 
The goal is to observe the impact of the mesh resolution on the amplitude growth.
Two dimensional simulations are 
performed on an uniform mesh with equal mesh spacing in $x$ and $y$ directions.   
The computational domain length is $L_x=8.9$~cm and $L_y=75$~cm.
The coarse, medium and fine grid simulations performed using the 
fifth order WENO method are on a uniform grid
that has $64, 128$ and $256$ grid points per initial perturbation wavelength.
Figure~\ref{fig:amp_sim3} shows the effect of the grid resolution on amplitude. 
The high-resolution simulation for $M=1.21$ is closest to the experimental measurements. 
The coarse grid simulation with
$64$ grid points per initial wavelength with a grid spacing
$\Delta x = \Delta y = 0.0928$~cm
leads to lower amplitudes. This is due to the low resolution in amplitude
where there are only $2$ grid points per initial amplitude.
This low resolution simulation does not capture the 
structure of the interface separating two fluids.
Figure~\ref{fig:fmc_densProfile} shows the need for
high-resolution simulations to accurately capture the fine-scale structure of the secondary
Kelvin-Helmholtz instability. 
The vortex roll-ups known as the mushroom-shaped structures are not observed in the coarse/medium simulations.
The key issue is the jump in the tangential component of the velocity across the interface.
This feature needs to be enabled to observe the secondary Kelvin-Helmholtz instability
and transport the vorticity along with the perturbed interface.
The simulations presented are performed at the Arkansas High Performance Computing Center where a standard node has Intel Xeon Gold 6130 CPUs with a total of 32 cores at a clock rate of 2.1 GHz.
Due to the limitation on the maximum node numbers per run,
the fine grid simulations are performed on 4 nodes with a total 128 processors.
To resolve the fine-scale structures, one can refine more and more the mesh. 
However, high resolution front tracking simulations with fifth-order WENO 
flux reconstructions and artificial compression of the diffuse interfaces 
become computationally very expensive. 

\section{Concluding Remarks}
\label{sec:conc}
Numerical simulations of turbulent mixing due to hydrodynamic
instabilities require an accurate representation of the perturbed interface between fluids. 
Accurate and robust front tracking simulations with the classical fifth order 
WENO scheme and artificial compression reveal agreement with experimental data.
The single-mode shock-induced RMI simulations of an air/SF$_6$ interface 
for the Mach numbers $M=1.11 $ and $M=1.21$ experiments of Collins and Jacobs 
show good agreement on the interface displacement 
and amplitude as a function of time. 
A difference less than 10\% is observed between the experiments and simulations 
interface velocities due to the expansion waves present only in the experiments.
Our validation studies show good agreement in early-time amplitude growth
between the fine grid simulation and the Mach number $M=1.11$ experiments.
However, we observe about 4\% discrepancy on the early-time growth rate for the Mach number $M=1.21$.  
For high Mach number,  
it is shown that the high resolution simulations with $256$ mesh points
per initial perturbation wavelength lie closer to the experimental measured data, 
however do not yet show numerical convergence at the grid resolutions presented in this paper.

In terms of repeatability of experiments, 
Collins and Jacobs reported differences (approximately 10\%) 
in amplitudes of the initial perturbations. 
According to Collins and Jacobs, 
increasing the interface amplitude results in a larger growth rate 
and more vorticity is created by the large interface amplitude values. 
The uncertainty quantification studies to investigate the effect of model
and input parameters on the growth rate and model improvement to capture 
the vortices at the interface are under development. 

\section*{Acknowledgments}
It is a great pleasure to thank Dr.~Jeffrey W.~Jacob for providing experimental data 
for detailed comparison studies, 
Dr.~James Glimm and Dr.~Xiaolin Li for many helpful comments and discussions.
Dr.Kaman acknowledges support of the Lawrence Jesser Toll Jr. Endowed Chair 
at the University of Arkansas.
This research is supported by the Arkansas High Performance Computing Center 
which is funded through multiple National Science Foundation grants and the 
Arkansas Economic Development Commission.


\end{document}